\documentstyle[amstex,amssymb]{article}

\newtheorem{fprop}{Proposition}[section]

\newtheorem{theor}{Theorem}[section]
\newtheorem{prop}{Proposition}[section]
\newtheorem{cor}{Corollary}[section]

\begin{document}

\author{Paul POPESCU}
\title{Higher order transverse bundles and riemannian foliations }
\date{}
\maketitle

\begin{abstract}
The purpose of this paper is to prove that each of the following conditions
is equivalent to that the foliation ${\cal F}$ is riemannian: 1) the lifted
foliation ${\cal F}^{r}$ on the $r$-transverse bundle $\nu ^{r}{\cal F}$ is
riemannian for an $r\geq 1$; 2) the foliation ${\cal F}_{0}^{r}$ on a
slashed $\nu _{\ast }^{r}{\cal F}$ is riemannian and vertically exact for an 
$r\geq 1$; 3) there is a positively admissible transverse lagrangian on a $%
\nu _{\ast }^{r}{\cal F}$, for an $r\geq 1$. Analogous results have been
proved previously for normal jet vector bundles.
\end{abstract}

\section{Introduction}

Most of geometrical objects considered on a differentiable manifold using
the tangent bundle can be constructed on a foliated manifold using the
normal bundle. This is described in an algebraic manner, using natural
functors, as in \cite{Wo3}. We use in this paper the normal bundle of order $%
r$, a foliated bundle that is a counterpart of the holonomy invariant of the
tangent bundle of order $r$.

Various conditions that a foliation be riemannian are studied in many
papers, for example \cite{JoWo, MiMo, MP09F, PMCR, Ta01, Wo2}.

The conditions studied in this paper are closely related to \cite{JoWo,
MiMo, MP09F, PMCR} and they have initially the origin in a special case of a
problem presented by E. Ghys in Appendix E of P. Molino's book \cite{Mo},
i.e. asking if the existence of a foliated Finsler metric assure that a
foliation is riemannian (see \cite{JoWo, MiMo, PMCR, MP09F} for more
details). According to \cite{MP09F}, the answer is affirmative in a more
general case of a transverse lagrangian fulfilling a natural regularity
condition, automatically fulfilled by a transverse Finslerian. The case of a
transverse Finslerian on a compact manifold is studied in \cite{JoWo}, using
a different method.

The conditions in \cite{MP09F, PMCR} involve the existence of suitable
admissible lagrangians or foliated metrics on the normal jet vector bundles.
The aim of this paper is to study similar problems as in \cite{PMCR}, but on
higher order normal bundles. We use basic notions about foliations from \cite%
{Mo, Wo3} and some notions relating foliated bundles and a basic result
stated in Proposition \ref{prop:23} from \cite{PMCR}. We use some basic
notions on higher order tangent spaces from \cite{Mi, Mo1, Mo2}; we point
some differences used here (for example, the lift of sections used here is
different from that used in \cite{Mo2} and the canonical inclusions in
Proposition \ref{prcanin}). Deeper properties of foliated hamiltonians are
studied in \cite{Va}; it can be a model in a further study of some geometric
objects considered in our paper.

The first goal of this paper is to find conditions that a foliation be
riemannian, involving general conditions on higher order normal bundles. But
some other aspects of the problem can be stressed. For example, the leaves
of a riemannian foliation ${\cal F}$ are compact, then the leaf spaces $M/%
{\cal F}$ is a Satake manifold (or a V-manifold, in the original terminology
of Satake), one of the first known non-trivial orbifold. The existence of a
transverse lagrangian or hamiltonian is worth to be studied on such
generalized manifolds, together with their physical properties; it is also
the case of the normal bundle of a foliation studied in \cite{MP09F}. By
Proposition \ref{prsemis} below, the existence of a foliated regular
lagrangian of order $r$, that has a positively defined hessian, gives rise
to a transverse riemannian metric on the normal bundle of order $r+1$ (see 
\cite{Mi} for the non-foliated case). Also, a transverse metric of a
riemannian foliation lifts to a transverse metric to the lifted foliation on
the normal bundle of order $r$, that becomes a riemannian foliation (see
Proposition \ref{pr0} for a simple construction). Thus it is natural to
consider the following problem: {\em under which conditions the existence of
a regular and positive lagrangian or of a transverse metric on the normal
bundle of a higher order transverse foliation assures that given foliation
is riemannian?} A similar problem is studied for normal jet vector bundles
(or $(p,r)$--velocities according to \cite{Wo2, Wo3}) in \cite{PMCR}.

The first main result proved in this paper asserts that {\em the lifted
foliation on a normal bundle of some order }$r\geq 1${\em \ is riemannian
iff the given foliation is riemannian} (Theorem \ref{thm:22}). A simple
consequence is that a foliation ${\cal F}$ is a riemannian one, provided
that the lifted foliation ${\cal F}^{r}$ on $\nu ^{r}{\cal F}$ is
transversaly parallelizable or almost parallelizable. The proof of Theorem %
\ref{thm:22} can not give any answer to the following question: {\em when is 
}${\cal F}${\em \ riemannian if the foliation induced on }$\nu ^{r}{\cal F}%
\backslash I_{r-1}^{r}(\nu ^{r-1}{\cal F})${\em \ is riemannian for some }$%
r\geq 1$? We give some answers to this question, as follows. {\em The lifted
foliation }$F_{0}^{r}${\em \ on a suitable slashed bundle }$\nu _{\ast
}^{r}F ${\em \ of the }$r${\em -normal bundle }$\nu ^{r}F${\em \ is
riemannian and vertically exact for some }$r\geq 1${\em \ iff }$F${\em \ is
riemannian} (Theorem \ref{thm:22-1}). In Theorem \ref{thm:21} we prove that 
{\em there is a positively admissible lagrangian on }$\nu ^{r}F${\em \ for
some }$r\geq 1${\em \ iff the foliation }$F${\em \ is riemannian}. These
three Theorems are analogous to \cite[Theorems 1.1, 1.2 and 1.3]{PMCR},
proved in the case of normal jet vector bundles.

As a conclusion, the results in this paper, together with that proved in 
\cite{PMCR} for normal jet vector bundles, confirm that asking some suitable
natural conditions on a higher order lagrangian, the given foliation is
necessarily riemannian; thus riemannian foliations are necessary setting to
study certain transverse lagrangians, subject to natural conditions,
considered on normal jet vector bundles or on higher order normal bundles of
a foliation.

\section{Basic notions and constructions related to the higher order
transverse bundles of a foliation}

Let $M$ be an $n$-dimensional manifold and ${\cal F}$ be a $k$-dimensional
foliation on $M$. We denote by $\tau {\cal F}$ and $\nu {\cal F}$ the
tangent plane field and the normal bundle respectively. A bundle $E$ over $M$
is called {\em foliated} if there is an atlas of local trivialisations on $E$
such that all the components of the structural functions are basic ones
(see, for example, \cite{Wo3}). In this case a canonical foliation ${\cal F}%
_{E}$ on $E$ is induced, having the same dimension $k$, such that $p$
restricted to leaves is a local diffeomorphism. In particular, we consider
affine and vector bundles that are foliated. For example, $\nu {\cal F}$ is
a foliated bundle and a natural foliation on $\nu {\cal F}$ can be
considered. According to \cite{PMCR}, a {\em positively admissible lagrangian%
} on a foliated vector bundle $p:E\rightarrow M$ is a continuous map $%
L:E\rightarrow I\!\!R$ that is asked to be differentiable at least when it
is restricted to the total space of the slashed bundle $E_{\ast
}=E\backslash \{\bar{0}\}\rightarrow M$, where $\{\bar{0}\}$ is the image of
the null section, such that the following conditions hold: 1) $L$ is
positively defined (i.e. its vertical hessian is positively defined) and $%
L(x,y)\geq 0=L(x,0)$, $(\forall )x\in M$ and $y\in E_{x}=p^{-1}(x)$; 2) $L$
is locally projectable on a transverse lagrangian; 3) there is a basic
function $\varphi :M\rightarrow (0,\infty )$, such that for every $x\in M$
there is $y\in E_{x}$ such that $L(x,y)=\varphi (x)$. If a positively
transverse lagrangian $F$ is $2$--homogeneous (i.e. $F(x,\lambda y)=\lambda
^{2}F(x,y)$, $(\forall )\lambda >0$), then $F$ is called a {\em finslerian};
it is also a positively admissible lagrangian, taking $\varphi \equiv 1$, or
any positive constant. For a foliated bundle, we can regard the vertical
bundle $VTE=\ker p_{\ast }\rightarrow E$ as a vector subbundle of $\nu {\cal %
F}_{E}\rightarrow E$ by mean of the canonical projection $TE\rightarrow \nu
F_{E}$, since $VTE$ is transverse to $\tau F_{E}$. Notice that if $%
p:E\rightarrow M$ is an affine bundle, then the vertical hessian ${\rm Hess}%
\ L$ of a lagrangian $L:E\rightarrow I\!\!R$ is a symmetric bilinear form on
the fibers of the vertical bundle $VTE$, given by the second order
derivatives of $L$, using the fiber coordinates (see \cite{MP09F} for more
details using coordinates).

In order to have a unitary form of the notions we use, we give now the
constructions of the higher order normal bundle of a foliation and the
related geometric objects used in the paper. They are similar to the
non-foliate (or foliated by points) case as in \cite{Mi, Mo1, Mo2}, some are
used for example in \cite{Wo1} in the foliate case.

Let us suppose that a foliation ${\cal F}$ is defined by a foliated atlas
having the generic coordinates $(x^{u},x^{\bar{u}})_{u=\overline{1,p},\bar{u}%
=\overline{1,q}}$ on a chart $(U,\varphi )$, where $p$ is the dimension of
the leaves and $q$ is the transverse dimension. The local form of the
pseudogroup according to the coordinates change, is $x^{u^{\prime
}}=x^{u^{\prime }}(x^{u},x^{\bar{u}})$, $x^{\bar{u}^{\prime }}=x^{\bar{u}%
^{\prime }}(x^{\bar{u}})$; here $(x^{u^{\prime }},x^{\bar{u}^{\prime }})$
are some coordinates in a chart $(U^{\prime },\varphi ^{\prime })$, $U\cap
U^{\prime }\neq \emptyset $, and the pseudogroup member is $\varphi ^{\prime
}\circ \varphi ^{-1}:\varphi (U\cap V)\rightarrow \varphi ^{\prime }(U\cap
V)\subset I\!\!R^{p}\times I\!\!R^{q}$. The coordinates on $TU$ are $%
(x^{u},x^{\bar{u}},y^{(1)u},y^{(1)\bar{u}})$ and the bundles $\tau {\cal F}%
_{U}$ and $\nu {\cal F}_{U}$ (restricted to $U)$, have as coordinates $%
(x^{u},x^{\bar{u}},y^{(1)u})$ and $(x^{u},x^{\bar{u}},y^{(1)\bar{u}})$
respectively.

Let us consider $x_{0}\in U$. Two curves $\gamma
_{1,2}:I_{a}=(-a,a)\rightarrow M$, $a>0$, have a {\em transverse contact of
order }$r\geq 0$ in $0$ if $\gamma _{1}(0)=\gamma _{2}(0)=x_{0}$ and if $r>0$%
, then $\frac{d^{j}(\gamma _{1}\circ \varphi ^{-1})^{\bar{u}}}{dt^{j}}(0)=%
\frac{d^{j}(\gamma _{2}\circ \varphi ^{-1})^{\bar{u}}}{dt^{j}}(0)$, $%
(\forall )\bar{u}=\overline{1,q}$, $j=\overline{1,r}$; we denote $\gamma _{1}%
\overset{r,x_{0}}{\symbol{126}{}}{}\gamma _{2}$. The ,,transverse contact of
order $r$\textquotedblright\ relation is an equivalence relation and all the
classes $\hat{\gamma}_{x_{0}}$, $x_{0}\in M$, give $\nu ^{r}{\cal F}$, i.e.
the {\em transverse space of order }$r$ of ${\cal F}$. The canonical
projection $\pi ^{r}:\nu ^{r}{\cal F}\rightarrow M$ gives a local trivial
bundle. The generic coordinates $(x^{u},x^{\bar{u}})$ on $U$ gives rise to
some generic coordinates $(x^{u},x^{\bar{u}},y^{(1)\bar{u}},\cdots ,y^{(r)%
\bar{u}})$ on $\nu ^{r}{\cal F}_{U}$, where $y^{(j)\bar{u}}=\frac{1}{r!}%
\frac{d^{j}(\gamma \circ \varphi ^{-1})^{\bar{u}}}{dt^{j}}(0)$ and $\hat{%
\gamma}_{x}$ is a class for $\overset{r,x}{\symbol{126}{}}{}$. Finally one
obtain a generic local chart $(\nu ^{r}{\cal F}_{U},\nu ^{r}\varphi )$ on $%
\nu ^{r}{\cal F}$ in a suitable atlas, according to every generic chart $%
(U,\varphi )$ on $M$, from a suitable atlas. We use notations and
constructions according to the non-foliate case, according to \cite[Sect. 6.1%
]{Mi} and \cite{MP06}.

The generic coordinates $(x^{u},x^{\bar{u}},y^{(1)\bar{u}},\cdots ,y^{(r)%
\bar{u}})$ on $\nu ^{r}{\cal F}_{U}$ change according to the rules: 
\begin{equation}
\begin{array}{clcr}
x^{u^{\prime }}= & x^{u^{\prime }}(x^{u},x^{\bar{u}}) &  &  \\ 
x^{\bar{u}^{\prime }}= & x^{\bar{u}^{\prime }}(x^{\bar{u}}) &  &  \\ 
y^{(1)\bar{u}^{\prime }}= & y^{(1)\bar{u}}{{{\displaystyle{\frac{\partial
x^{i^{\prime }}}{\partial x^{i}}}}}\qquad } &  &  \\ 
2y^{(2)\bar{u}^{\prime }}= & y^{(1)\bar{u}}{{{\displaystyle{\frac{\partial
y^{(1)\bar{u}^{\prime }}}{\partial x^{\bar{u}}}}}}}+ & 2y^{(2)\bar{u}}{{{%
\displaystyle{\frac{\partial y^{(1)\bar{u}^{\prime }}}{\partial y^{(1)\bar{u}%
}}}}}}, &  \\ 
\vdots &  &  &  \\ 
ry^{(r)\bar{u}^{\prime }}= & y^{(1)\bar{u}}{{{\displaystyle{\frac{\partial
y^{(r-1)\bar{u}^{\prime }}}{\partial x^{\bar{u}}}}}}}+ & \cdots & +ry^{(r)%
\bar{u}}{{{\displaystyle{\frac{\partial y^{(r-1)\bar{u}^{\prime }}}{\partial
y^{(r-1)\bar{u}}}}}}}.%
\end{array}
\label{0eqchy}
\end{equation}

Denoting $x^{\bar{u}}=y^{(0)\bar{u}}$, one have, for $0\leq \alpha \leq
\beta \leq \gamma \leq r$:%
\[
\frac{\partial y^{(\gamma )\bar{u}^{\prime }}}{\partial y^{(\beta )\bar{u}}}=%
\frac{\partial y^{(\gamma -\alpha )\bar{u}^{\prime }}}{\partial y^{(\beta
-\alpha )\bar{u}}}. 
\]

In particular,%
\[
{{{{\frac{\partial y^{(r-1)\bar{u}^{\prime }}}{\partial y^{(r-1)\bar{u}}}=}}}%
}\dfrac{\partial x^{\bar{u}^{\prime }}}{\partial x^{\bar{u}}}. 
\]

These formulas are true according to similar ones in the non-foliate case 
\cite[Sect. 6.1]{Mi}.

Various bundle structures can be considered over a $\nu ^{r}{\cal F}$; for
example, for $0\leq r^{\prime }\leq r$, the canonical projection $\pi
_{r^{\prime }}^{r}:\nu ^{r}{\cal F}\rightarrow \nu ^{r^{\prime }}{\cal F}$
is a foliated bundle. In particular, for $r\geq 1$, $\pi _{r-1}^{r}:\nu ^{r}%
{\cal F}\rightarrow \nu ^{r-1}{\cal F}$ is a (foliated) affine bundle for $%
r>1$ and $\pi _{0}^{1}:\nu {\cal F}\rightarrow \nu ^{0}{\cal F}=M$ is a
(foliated) vector bundle (for $r=1$).

\begin{fprop}
\label{prcanin}For $1\leq r^{\prime }\leq r$, there is an inclusions of
foliated submanifolds (in fact of foliated subbundles over $M$), $%
I_{r^{\prime }}^{r}:\nu ^{r^{\prime }}{\cal F}\rightarrow \nu ^{r}{\cal F}$,
where the inclusion assigns to an equivalence class in $[\gamma ]\in \nu
_{;m}^{r^{\prime }}{\cal F}$ an equivalence class in $\nu _{;m}^{r}{\cal F}$
that the first $r-r^{\prime }$ derivatives vanish, then the next $r^{\prime
} $ derivatives are the same as the first $r^{\prime }$ derivatives of $%
\gamma $.
\end{fprop}

{\em Proof.} We use generic coordinates. Denoting $z^{(\alpha )\bar{u}%
}=\alpha !y^{(\alpha )\bar{u}}$, $\alpha =\overline{1,r}$ as new
coordinates, the local form of $I_{r^{\prime }}^{r}$ is 
\[
(x^{u},x^{\bar{u}},z^{(1)\bar{u}},\cdots ,z^{(r^{\prime })\bar{u}})\overset{%
I_{r^{\prime }}^{r}}{\rightarrow }{}(x^{u},x^{\bar{u}},0,\ldots ,0,z^{(1)%
\bar{u}},\cdots ,z^{(r^{\prime })\bar{u}}); 
\]
using again generic coordinates, the local form of $I_{r^{\prime }}^{r}$ is 
\[
(x^{u},x^{\bar{u}},y^{(1)\bar{u}},\cdots ,y^{(r^{\prime })\bar{u}})\overset{%
I_{r^{\prime }}^{r}}{\rightarrow }{}(x^{u},x^{\bar{u}},0,\ldots ,0,\frac{%
(r-r^{\prime }+1)!}{1!}y^{(1)\bar{u}},\cdots ,\frac{r!}{r^{\prime }!}%
y^{(r^{\prime })\bar{u}}). 
\]
Using formulas (\ref{0eqchy}), it easily follows that $I_{r^{\prime }}^{r}$
is globally defined. $\Box $

Let us notice that the local form of $I_{r-1}^{r}$ is 
\[
(x^{u},x^{\bar{u}},y^{(1)\bar{u}},\cdots ,y^{(r-1)\bar{u}})\overset{%
I_{r^{\prime }}^{r}}{\rightarrow }{}(x^{u},x^{\bar{u}},0,\frac{2!}{1!}y^{(1)%
\bar{u}},\cdots ,\frac{r!}{(r-1)!}y^{(r^{\prime })\bar{u}}) 
\]
and $I_{r}^{r}$ is the identity of $\nu ^{r}{\cal F}$.

Thus we have $I_{0}^{r}(M)\subset I_{1}^{r}(\nu {\cal F})\subset
I_{2}^{r}(\nu ^{2}{\cal F})\subset \cdots \subset I_{r-1}^{r}(\nu ^{r-1}%
{\cal F})\subset \nu ^{r}{\cal F}$.

A transverse vector field $\bar{X}\in \Gamma (\nu {\cal F})$ lifts in this
way to the transverse section $I_{1}^{r}(\bar{X}):M\rightarrow \nu ^{r}{\cal %
F}$ of the bundle $\pi _{0}^{r}:\nu ^{r}{\cal F}\rightarrow M$. A more
simple case is when $\bar{X}=\bar{0}$ is the null vector field; its lift is
the section $\bar{0}^{r}:M\rightarrow \nu ^{r}{\cal F}$, $\bar{0}%
^{r}(m)=I_{0}^{r}(m)$. Notice that, in the non-foliate case, the lifts
constructed here are not the same as the lifts constructed in \cite{Mo2},
where the lifts are vector fields.

We denote by ${\cal F}^{r}$ the foliation on $\nu ^{r}{\cal F}$. In a
similar way as in the non-foliate case in \cite[Sect. 6.1]{Mi}, we perform
some constructions that are useful later. Notice that in the foliate case
the transverse $\nu {\cal F}^{r}$ play the role of a tangent space for $\nu
^{r}{\cal F}$, as the tangent space $\tau T^{r}M$ is for $T^{r}M$ in the
non-foliate case in \cite{Mi}.

For every $r\geq 1$ and $0\leq r^{\prime }\leq r$, , the canonical
projection $\pi _{r^{\prime }}^{r}:\nu ^{r}{\cal F}\rightarrow \nu
^{r^{\prime }}{\cal F}$ induces a transverse map $\bar{\pi}_{r^{\prime
}}^{r}:\nu {\cal F}^{r}\rightarrow \nu {\cal F}^{r^{\prime }}$ that is a
vector bundle map of foliated vector bundles; notice that $\pi _{0}^{r}=\pi
^{r}$, ${\cal F}^{0}={\cal F}$, $\nu ^{1}{\cal F}=\nu {\cal F}$, $\nu ^{0}%
{\cal F}=M$ and $\bar{\pi}_{0}^{r}=\bar{\pi}^{r}$. We denote the kernel
vector subbundle bundle $\ker \bar{\pi}_{r^{\prime }}^{r}\subset \nu {\cal F}%
^{r}$ by $\bar{V}_{r^{\prime }}^{r}$; it is a foliate vector bundle as well.
Since for $r_{1}\leq r_{2}\leq r_{3}$, one have $\pi _{r_{1}}^{r_{3}}=\pi
_{r_{2}}^{r_{3}}\circ \pi _{r_{1}}^{r_{2}}$ and $\bar{\pi}_{r_{1}}^{r_{3}}=%
\bar{\pi}_{r_{2}}^{r_{3}}\circ \bar{\pi}_{r_{1}}^{r_{2}}$, it follows that
there are foliated vector subbundles $\bar{V}_{r-1}^{r}\subset \bar{V}%
_{r-2}^{r}\subset \cdots \subset \bar{V}_{0}^{r}\subset \nu {\cal F}^{r}$.
Notice that $\nu ^{r+1}{\cal F}\subset \nu {\cal F}^{r}$ is an affine
subbundle over $\nu ^{r}{\cal F}$, for $r\geq 1$, while $\nu ^{1}{\cal F}%
=\nu {\cal F}^{0}=\nu {\cal F}$ for $r=0$.

There is an $r$-transverse structures in the fibers of on $\nu {\cal F}^{r}$%
, i.e. a vector bundle map $J:\nu {\cal F}^{r}\rightarrow \nu {\cal F}^{r}$
(analogous of the $r$-tangent structures in the non-foliate case), and its
dual $J^{\ast }:\nu ^{\ast }{\cal F}^{r}\rightarrow \nu ^{\ast }{\cal F}^{r}$%
, given in generic local coordinates by%
\begin{eqnarray*}
J &=&\overline{\dfrac{\partial }{\partial y^{(1)\bar{u}}}}\otimes \overline{%
dx^{\bar{u}}}+\overline{\dfrac{\partial }{\partial y^{(2)\bar{u}}}}\otimes 
\overline{dy^{(1)\bar{u}}}+\cdots +\overline{\dfrac{\partial }{\partial
y^{(r)\bar{u}}}}\otimes \overline{dy^{(r-1)\bar{u}}}, \\
J^{\ast } &=&\overline{dx^{\bar{u}}}\otimes \overline{\dfrac{\partial }{%
\partial y^{(1)\bar{u}}}}+\overline{dy^{(1)\bar{u}}}\otimes \overline{\dfrac{%
\partial }{\partial y^{(2)\bar{u}}}}+\cdots +\overline{dy^{(r-1)\bar{u}}}%
\otimes \overline{\dfrac{\partial }{\partial y^{(r)\bar{u}}}}.
\end{eqnarray*}

A {\em transverse }$r${\em -non-linear connection} is a splitting of $\nu 
{\cal F}^{r}$ as a Whitney sum of transverse vector bundles%
\begin{equation}
\nu {\cal F}^{r}=\bar{V}_{0}^{r}\oplus \bar{H}_{0}^{r},  \label{eqtran01}
\end{equation}%
where $\bar{H}_{0}^{r}$ is the $r$-horizontal vector bundle, that is
canonically isomorphic with $(\bar{\pi}^{r})^{\ast }\nu {\cal F}$. We denote
by $h:\nu {\cal F}^{r}\rightarrow \bar{H}_{0}^{r}$ the projector given by
the above decomposition. Using generic local coordinates, the local form (on
the fibers) of the projector $h$ is 
\begin{equation}
(X^{\bar{u}},Y^{(1)\bar{u}},\cdots ,Y^{(r)\bar{u}})\overset{h}{\rightarrow }%
{}(X^{\bar{u}}+N_{(1)\bar{v}}^{\bar{u}}Y^{(1)\bar{v}}+\cdots +N_{(r)\bar{v}%
}^{\bar{u}}Y^{(r)\bar{v}}).  \label{eqcoordconn}
\end{equation}

Given a transverse $r$-non-linear connection by a splitting (\ref{eqtran01}%
), the consecutive images by $J$ in the fibers of $\nu {\cal F}^{r}$,%
\[
J\left( \bar{H}_{0}^{r}\right) =\bar{H}_{1}^{r},\ldots ,J\left( \bar{H}%
_{r-1}^{r}\right) =\bar{H}_{r}^{r} 
\]%
define some transverse vector subbundles of $\nu {\cal F}^{r}$, all
isomorphic with $\bar{H}_{0}^{r}$, such that there are the following Whitney
sum decompositions 
\begin{equation}
\bar{V}_{0}^{r}=\bar{H}_{1}^{r}\oplus \cdots \oplus \bar{H}_{r}^{r},\;\nu 
{\cal F}^{r}=\bar{H}_{0}^{r}\oplus \bar{H}_{1}^{r}\oplus \cdots \oplus \bar{H%
}_{r}^{r}.  \label{eqtran01a}
\end{equation}

Notice that $\bar{H}_{r}^{r}=\bar{V}_{r-1}^{r}$ and we can prove the
following result.

\begin{prop}
\label{prsplit01}Any splitting $\nu {\cal F}^{r}=\bar{V}_{r-1}^{r}\oplus 
\bar{H}_{r-1}^{r}$ gives rise to a splitting (\ref{eqtran01}).
\end{prop}

{\em Proof.} The given splitting gives rise also to a splitting $\nu ^{\ast }%
{\cal F}^{r}=\left( \bar{V}_{r-1}^{r}\right) ^{\ast }\oplus \left( \bar{H}%
_{r-1}^{r}\right) ^{\ast }$ and denote $\bar{W}^{r}\overset{not.}{=}{}\left( 
\bar{V}_{r-1}^{r}\right) ^{\ast }\subset \nu ^{\ast }{\cal F}^{r}$. Let us
consider the consecutive images by $J^{\ast }$ in the fibers of $\nu ^{\ast }%
{\cal F}^{r}$,%
\[
J^{\ast }\left( \bar{W}^{r}\right) =\bar{W}^{r-1},\ldots ,J^{\ast }\left( 
\bar{W}^{1}\right) =\bar{W}^{0} 
\]%
defining some transverse vector subbundles of $\nu ^{\ast }{\cal F}^{r}$,
all isomorphic with $\left( \bar{V}_{r-1}^{r}\right) ^{\ast }$ and where $%
\bar{W}^{0}=\left( \bar{V}_{0}^{r}\right) ^{\circ }$ is the polar (or the
annihilator). Considering $\bar{W}=\bar{W}^{1}\oplus \cdots \oplus \bar{W}%
^{r}$, we have $\nu ^{\ast }{\cal F}^{r}=\left( \bar{V}_{0}^{r}\right)
^{\circ }\oplus \bar{W}$, thus $\bar{W}$ is isomorphic with $\left( \bar{V}%
_{0}^{r}\right) ^{\ast }$. It follows that the inclusion $\left( \bar{V}%
_{0}^{r}\right) ^{\ast }\subset \nu ^{\ast }{\cal F}^{r}$ reverses by
duality to a transverse epimorphism (in fact a projector) $\Pi :\nu {\cal F}%
^{r}\rightarrow \bar{V}_{0}^{r}$ that gives a splitting (\ref{eqtran01}),
where $\bar{H}_{0}^{r}=\ker \Pi $. $\Box $

Using generic local coordinates, the local form (on the fibers) of the
projector $\Pi :\nu {\cal F}^{r}\rightarrow \bar{V}_{r-1}^{r}$ is 
\begin{equation}
(X^{\bar{u}},Y^{(1)\bar{u}},\cdots ,Y^{(r)\bar{u}})\overset{\Pi }{%
\rightarrow }{}(M_{(0)\bar{v}}^{\bar{u}}X^{\bar{v}}+M_{(1)\bar{v}}^{\bar{u}%
}Y^{(1)\bar{v}}\cdots +M_{(r-1)\bar{v}}^{\bar{u}}Y^{(r-1)\bar{v}}+Y^{(r)\bar{%
u}}).  \label{eqcoordconn02}
\end{equation}

The local coefficients $N$ and $M$ in formulas (\ref{eqcoordconn}) and (\ref%
{eqcoordconn02}) respectively are called as {\em dual coefficients} in \cite[%
sect. 6.6]{Mi}.

A {\em transverse }$r${\em -semi-spray} is a foliate section $S:\nu ^{r}%
{\cal F}\rightarrow \nu ^{r+1}{\cal F}$ of the affine bundle $\pi
_{r}^{r+1}:\nu ^{r+1}{\cal F}\rightarrow \nu ^{r}{\cal F}$. Since $\nu ^{r+1}%
{\cal F}\subset \nu {\cal F}^{r}$, it follows that an $r$-semi-spray can be
regarded as well as a transverse section $S:\nu ^{r}{\cal F}\rightarrow \nu 
{\cal F}^{r}$; using generic local coordinates, the local form of $S:\nu ^{r}%
{\cal F}\rightarrow \nu ^{r+1}{\cal F}$ is 
\begin{equation}
(x^{u},x^{\bar{u}},y^{(1)\bar{u}},\cdots ,y^{(r)\bar{u}})\overset{S}{%
\rightarrow }{}(x^{u},x^{\bar{u}},y^{(1)\bar{u}},\cdots ,ry^{(r)\bar{u}%
},(r+1)S^{\bar{u}}(x^{\bar{u}},y^{(1)\bar{u}},\cdots ,y^{(r)\bar{u}})).
\label{eqlocS}
\end{equation}

\begin{prop}
\label{prnlc}Any transverse $r$-semi-spray gives rise to a transvese $r$%
-non-linear connection, i.e. a splitting (\ref{eqtran01}).
\end{prop}

{\em Proof.} If $\bar{X}\in \Gamma (\nu {\cal F}^{r})$ is a transverse
section, let us denote by $L_{\bar{X}}$ the transverse Lie derivation
induced in the transverse tensor algebra. As in the non-foliate case in \cite%
{CSC1}, the operators 
\begin{equation}
h=\dfrac{1}{k+1}\left( kI-L_{S}J\right) ,v=\dfrac{1}{k+1}\left(
I+L_{S}J\right)  \label{invar01}
\end{equation}%
are complementary projectors on the fibers of the transverse bundle $\nu 
{\cal F}^{r}\rightarrow \nu ^{r}{\cal F}$ and $v(\nu {\cal F}^{r})=\bar{V}%
_{r-1}^{r}$. Using Proposition \ref{prsplit01}, a transvese $r$-non-linear
connection\ is obtained. $\Box $

Using generic coordinates, it can be proved using (\ref{invar01}), that the
dual coefficients (\ref{eqcoordconn02}) of the non-linear connection are
given by the formulas%
\begin{equation}
M_{(r)\bar{v}}^{\bar{u}}=-\dfrac{\partial S^{\bar{u}}}{\partial y^{(1)\bar{v}%
}},M_{(r-1)\bar{v}}^{\bar{u}}=-\dfrac{\partial S^{\bar{u}}}{\partial y^{(2)%
\bar{v}}},M_{(1)\bar{v}}^{\bar{u}}=-\dfrac{\partial S^{\bar{u}}}{\partial
y^{(r)\bar{v}}}.  \label{eqdefsemis}
\end{equation}

A fact that we use latter is the following result.

\begin{prop}
\label{prtrmet01}A transvese $r$-non-linear connection{\em \ }and a
transverse riemannian metric in the fibers of $\bar{V}_{r-1}^{r}$ lifts to a
transverse riemannian metric on $\nu {\cal F}^{r}$. Conversely, a transverse
riemannian metric on $\nu {\cal F}^{r}$ gives a transvese $r$-non-linear
connection{\em \ }and a transverse riemannian metric in the fibers of $\bar{V%
}_{r-1}^{r}$.
\end{prop}

{\em Proof.} Every riemannian metric in the fibers of $\bar{V}%
_{r-1}^{r}\cong \bar{H}_{0}^{r}\cong \left( \bar{\pi}_{0}^{r}\right) ^{\ast
}\nu {\cal F}$ lifts to transverse riemannian metrics on $\bar{H}_{1}^{r}$, $%
\ldots ,\bar{H}_{r}^{r}$ and consequently to a riemannian metric on $\nu 
{\cal F}^{r}$, that becomes a riemannian foliation. Conversely, if ${\cal F}%
^{r}$ is a riemannian foliation, then a transverse riemannian metric in the
fibers of $\nu {\cal F}^{r}$ gives a decomposition (\ref{eqtran01}), where $%
\bar{H}_{0}^{r}=\left( \bar{V}_{0}^{r}\right) ^{\bot }$ and induces a
transverse riemannian metric in the fibers of the vector subbundle $\bar{H}%
_{0}^{r}$. $\Box $

The $r$-transverse non-linear connections, semi-sprays and riemannian
metrics are involved in the case of regular $r$-transverse lagrangians that
we consider in the sequel.

An $r${\em -transverse lagrangian} (a transverse lagrangian of order $r\geq
1 $, i.e. locally projectable on an $r$-lagrangian) is a continuous real map 
$L:\nu ^{r}{\cal F}\rightarrow I\!\!R,$ smooth on an open fibered
submanifold $\nu _{\ast }^{r}{\cal F\subset \nu }^{r}{\cal F}$. The cases
studied in the paper are when $\nu _{\ast }^{r}{\cal F=}\nu ^{r}{\cal F}$,
i.e. $L$ is smooth, or when $\nu ^{r}{\cal F}\backslash \nu _{\ast }^{r}%
{\cal F}$ contains $I_{r-1}^{r}(\nu ^{r-1}{\cal F})$, i.e. $L$ is slashed.
For sake of simplicity, we perform the next constructions in the case of a
smooth $L$, in the slashed case we must be care of domains where the objects
are defined. The {\em vertical hessian} of $L$ is the bilinear form $h$ in
the fibers of $\bar{V}_{r-1}^{r}$, given in some generic coordinates by%
\[
h_{\bar{u}\bar{v}}=\dfrac{\partial ^{2}L}{\partial y^{(r)\bar{u}}\partial
y^{(r)\bar{v}}}. 
\]%
We say that $L$ is {\em regular} if its vertical hessian is non-degenerated.
The fibers of the fibered manifold $\nu ^{r}{\cal F}\rightarrow \nu ^{r-1}%
{\cal F}$ are affine spaces. Using (\ref{0eqchy}), the generic coordinates
on fibers change according the formulas%
\[
ry^{(r)\bar{u}^{\prime }}=\Gamma (y^{(r-1)\bar{u}^{\prime }})+ry^{(r)\bar{u}}%
\dfrac{\partial x^{\bar{u}^{\prime }}}{\partial x^{\bar{u}}}{,} 
\]%
where%
\[
\Gamma =y^{(1)\bar{u}}\dfrac{\partial }{\partial x^{\bar{u}}}+2y^{(2)\bar{u}}%
\dfrac{\partial }{\partial y^{(1)\bar{u}}}+\cdots +ry^{(r)\bar{u}}\dfrac{%
\partial }{\partial y^{(r-1)\bar{u}}}. 
\]

\begin{prop}
\label{prsemis}a) If an $r$-lagrangian $L$ is regular, then it can define
canonically a transverse $r$-semi-spray and a transverse $r$-non- linear
connection.

b) If the vertical hessian of an $r$-lagrangian $L$ is positively defined,
then ${\cal F}^{r}$ is a riemannian foliation.
\end{prop}

{\em Proof.} In order to prove a) it suffices to construct a transverse $r$%
-semi-spray. Then a transverse $r$-non- linear connection can be constructed
by Proposition \ref{prnlc}. We use generic local coordinates, according to
formula (\ref{eqlocS}). As in the non-foliate case \cite[Theorem 8.8.1]{Mi}
(see also \cite{Bu}), the local form of the functions $S^{\bar{u}}$ can be
taken according to the formula%
\[
S^{\bar{u}}=\frac{1}{2\left( r+1\right) }h^{\bar{u}\bar{v}}\left( \Gamma
\left( \dfrac{\partial L}{\partial y^{(r)\bar{v}}}\right) -\dfrac{\partial L%
}{\partial y^{(r-1)\bar{v}}}\right) , 
\]%
where $\left( h^{\bar{u}\bar{v}}\right) =\left( h_{\bar{u}\bar{v}}\right)
^{-1}$. In order to prove b), one construct a transverse riemannian metric $%
H $ in the fibers of $\nu {\cal F}^{r}$, using its decomposition given by
formula (\ref{eqtran01a}) and taking into account that all the summands are
isomorphic with $V_{r-1}^{r}$, where $h$ is a riemannian metric on fibers. $%
\Box $

According to the case of trivial foliation of $M$ by points in \cite{MP06}, $%
\nu ^{r-1}{\cal F}\times _{M}\nu ^{\ast }{\cal F}\overset{not.}{=}{}\nu
^{r\ast }{\cal F}$ play the role of the vectorial dual of the affine bundle $%
\nu ^{r}{\cal F}\rightarrow \nu ^{r-1}{\cal F}$. The usual partial
derivatives of $L$ in the highest order transverse coordinates define a
well-defined Legendre map ${\cal L}:\nu ^{r}\rightarrow \nu ^{r\ast }{\cal F}
$, i.e.%
\[
(x^{u},x^{\bar{u}},y^{(1)\bar{u}},\cdots ,y^{(r)\bar{u}})\overset{{\cal L}}{%
\rightarrow }{}(x^{u},x^{\bar{u}},y^{(1)\bar{u}},\cdots ,y^{(r-1)\bar{u}},%
\frac{\partial L}{\partial y^{(r)\bar{u}}}). 
\]

If $L$ is regular, then ${\cal L}$ is a local diffeomorphism. If ${\cal L}$
is a global diffeomorphism we say that $L$ is {\em hyperregular}. We say
that $H:\nu ^{r\ast }{\cal F}\rightarrow I\!\!R$, $H=L\circ {\cal L}^{-1}$
is the pseudo-hamiltonian associated with $L$. If ${\cal L}^{-1}$ has the
local form%
\[
{\cal L}^{-1}(x^{u},x^{\bar{u}},y^{(1)\bar{u}},\cdots ,y^{(r-1)\bar{u}},p_{%
\bar{u}})=(x^{u},x^{\bar{u}},y^{(1)\bar{u}},\cdots ,y^{(r-1)\bar{u}},H^{\bar{%
u}}(x^{u},x^{\bar{u}},y^{(1)\bar{u}},\cdots ,y^{(r-1)\bar{u}},p_{\bar{u}})), 
\]%
then $H$ has the local form%
\[
H(x^{u},x^{\bar{u}},y^{(1)\bar{u}},\cdots ,y^{(r-1)\bar{u}},p_{\bar{u}%
})=L(x^{u},x^{\bar{u}},y^{(1)\bar{u}},\cdots ,y^{(r-1)\bar{u}},H^{\bar{u}%
}(x^{u},x^{\bar{u}},y^{(1)\bar{u}},\cdots ,y^{(r-1)\bar{u}},p_{\bar{u}})). 
\]

For $0\leq r^{\prime }\leq r$, let us denote $\nu ^{r^{\prime },(r-r^{\prime
})\ast }{\cal F}=\nu ^{r^{\prime }}{\cal F}\times _{M}\left( \nu ^{\ast }%
{\cal F}\right) ^{r-r^{\prime }}$, where $\left( \nu ^{\ast }{\cal F}\right)
^{r-r^{\prime }}=$ $\nu ^{\ast }{\cal F}\times _{M}\cdots \times _{M}\nu
^{\ast }{\cal F}$, with the fibered product of$\ (r-r^{\prime })$-times. In
particular, $\nu ^{r\ast }=\nu ^{r-1,r\ast }{\cal F}=\nu ^{r-1}{\cal F}%
\times _{M}\nu ^{\ast }{\cal F}$.

A {\em transverse slashed lagrangian} of order $r$ is a continuous map $%
L^{r}:$ $\nu ^{r}{\cal F}\rightarrow I\!\!R$ that is differentiable on an
open fibered submanifold $\nu _{\ast }^{r}{\cal F\subset }\nu ^{r}{\cal F}$,
called a {\em slashed bundle}. All the above constructions can be adapted
for slashed lagrangians.

Let us suppose that $L^{r}$ is {\em hyperregular}, i.e. the Legendre map $%
{\cal L}^{(r)}:\nu _{\ast }^{r}\rightarrow \nu ^{1,(r-1)\ast }{\cal F}=\nu
^{r-1}{\cal F}\times _{M}\nu ^{\ast }{\cal F}$ is a diffeomorphism on its
image Let us suppose also that ${\cal L}^{(r)}\left( \nu _{\ast }^{r}\right)
=$ $\nu _{\ast }^{1,(r-1)\ast }{\cal F}=$ $\nu _{\ast }^{r-1}{\cal F}\times
_{M}\nu _{\ast }^{\ast }{\cal F}$; here $\nu _{\ast }^{\ast }{\cal F}=\nu
^{\ast }{\cal F}\backslash \{\bar{0}\}$ (where $\{\bar{0}\}$ is the image of
the section obtained by all velocities vanishing, can be identified with $M$%
) and $\nu _{\ast }^{r-1}{\cal F}$ is a slashed subbundle of $\nu ^{r-1}%
{\cal F}$. We denote by $H^{1,r-1}=L^{r}\circ \left( {\cal L}^{(r)}\right)
^{-1}:\nu _{\ast }^{1,(r-1)\ast }{\cal F}\rightarrow I\!\!R$ its
pseudo-hamiltonian. (See \cite{MP06} for its classical definition and \cite%
{MP09B} for a coordinate description of the whole construction in the
non-foliate case). Analogous, for $0\leq j<r-1$, we suppose, step by step,
backward from $r-1$ from $0$, that there the usual partial derivatives of $%
L^{(j+1)}:\nu _{\ast }^{j+1,(r-j-1)\ast }{\cal F}=\nu _{\ast }^{r-j-1}{\cal F%
}\times _{M}\left( \nu _{\ast }^{\ast }{\cal F}\right) ^{j+1}\rightarrow
I\!\!R$ in the highest order transverse coordinates (of order $j+1$) define
a well-defined Legendre map ${\cal L}^{(j+1)}:\nu _{\ast }^{j+1,(r-j-1)\ast }%
{\cal F}=\nu _{\ast }^{j+1}{\cal F}\times _{M}\left( \nu _{\ast }^{\ast }%
{\cal F}\right) ^{r-j-1}\rightarrow \nu ^{j,(r-j)\ast }{\cal F}=\nu ^{j}%
{\cal F}\times _{M}\left( \nu ^{\ast }{\cal F}\right) ^{r-j}$. We suppose
that ${\cal L}^{(j+1)}$ is a diffeomorphism on its image and the image is
exactly ${\cal L}^{(j+1)}\left( \nu _{\ast }^{j+1,(r-j-1)\ast }{\cal F}%
\right) =\nu _{\ast }^{j,(r-j)\ast }{\cal F}=\nu _{\ast }^{j}{\cal F}\times
_{M}\left( \nu _{\ast }^{\ast }{\cal F}\right) ^{r-j}$. Then the
pseudo-hamiltonian $L^{(j)}=L^{(j+1)}\circ \left( {\cal L}^{(j+1)}\right)
^{-1}:\nu _{\ast }^{j,(r-j)\ast }{\cal F}\rightarrow I\!\!R$ can be
considered. Finally, for $j=0$, we obtain a transverse slashed lagrangian $%
L^{(0)}=L^{1}\circ \left( {\cal L}^{(1)}\right) ^{-1}:\nu _{\ast }^{0,r\ast }%
{\cal F}=\left( \nu _{\ast }^{\ast }{\cal F}\right) ^{r}\rightarrow I\!\!R$
and we suppose that ${\cal L}^{(1)}:\nu _{\ast }^{1,(r-1)\ast }{\cal F}=$ $%
\nu _{\ast }{\cal F}\times _{M}\left( \nu _{\ast }^{\ast }{\cal F}\right)
^{r-1}\rightarrow \nu _{\ast }^{0,r\ast }{\cal F}=$ $\left( \nu _{\ast
}^{\ast }{\cal F}\right) ^{r}\subset $ $\nu ^{0,r\ast }{\cal F}=$ $\left(
\nu ^{\ast }{\cal F}\right) ^{r}$ is a diffeomorphism. It follows a
diffeomorphism ${\cal L=L}^{(1)}\circ \cdots \circ {\cal L}^{(r)}:\nu _{\ast
}^{r}\rightarrow \left( \nu _{\ast }^{\ast }{\cal F}\right) ^{r}$ and a
transverse slashed lagrangian $L^{(0)}:\left( \nu _{\ast }^{\ast }{\cal F}%
\right) ^{r}\rightarrow I\!\!R$. The canonical diagonal inclusion $\nu
^{\ast }{\cal F}\rightarrow \left( \nu ^{\ast }{\cal F}\right) ^{r}$ sends $%
\nu _{\ast }^{\ast }{\cal F}\rightarrow \left( \nu _{\ast }^{\ast }{\cal F}%
\right) ^{r}$. We suppose that the restriction of $L^{(0)}$ to the diagonal
is a positively admissible lagrangian on $\nu ^{\ast }{\cal F}$, in fact a
transverse hamiltonian $H:\nu _{\ast }^{\ast }{\cal F}\rightarrow I\!\!R$.
If the given transverse lagrangian $L^{r}:\nu ^{r}{\cal F}\rightarrow I\!\!R$
fulfills all the above conditions, we say that $L$ itself is a {\em %
positively admissible lagrangian} (of order $r$) and $H$ is its{\em \
diagonal} {\em hamiltonian}.

The existence of a lifted metric, from the base space to the higher order
tangent bundle, is an well-known fact in the non-foliate case (see, for
example \cite[Sect. 9.2]{Mi}); we have to consider a simpler construction in
the foliated case, that it is also vertically exact, as in \cite{MP06, MP09B}%
.

\begin{prop}
\label{pr0}Any transverse metric $g$ on $\nu F$ gives canonically a
positively admissible lagrangian $L^{(r)}$ of order $r$ and a canonical
vertically exact transverse riemannian metric $g^{(r)}$ on $\nu ^{r}{\cal F}$%
, for any $r\geq 1$.
\end{prop}

{\em Proof.} We proceed by induction over $r\geq 1$. We consider the
quadratic first order lagrangian $L^{(1)}:\nu {\cal F}\rightarrow I\!\!R$, $%
L^{(1)}(x,y^{(1)})=g_{x}(y^{(1)},y^{(1)})$. The Levi-Civita connection of
the transverse metric $g$ on $\nu F$ gives rise to the geodesic first order
spray $S^{(1)}:\nu {\cal F}\rightarrow \nu ^{2}{\cal F}$ of $L^{(1)}$ and to
a second order lagrangian $L^{(2)}:\nu ^{2}{\cal F}\rightarrow I\!\!R$, $%
L^{(2)}(x,y^{(1)},y^{(2)})=L^{(1)}(x,y^{(1)})+L^{(1)}(x,y^{(2)}-S^{(1)}(x,y^{(1)})) 
$. Assume that $L^{(r-1)}:\nu ^{r-1}{\cal F}\rightarrow I\!\!R$ has been
constructed. An $(r-1)$-order spray $S^{(r)}:\nu ^{r}{\cal F}\rightarrow \nu
^{r+1}{\cal F}$ can be constructed according to Proposition \ref{prsemis},
since $L^{(r-1)}$ is $r$-regular. It follows $L^{(r)}:\nu ^{r}{\cal F}%
\rightarrow I\!\!R$, $L^{(r)}(x,y^{(1)},\ldots
,y^{(r)})=L^{(r-1)}(x,y^{(1)},\ldots
,y^{(r-1)})+L^{(1)}(x,y^{(r)}-S^{(r-1)}) $, that is a positively admissible
lagrangian of order $r$ ; the diagonal hamiltonian of $L^{(r)}$ is just the
dual hamiltonian of $L^{(1)}$. According to Proposition \ref{prtrmet01}, the
lagrangian $L^{(r)}$ gives rise to a transverse riemannian metric in the
fibers of $\nu {\cal F}^{r}$, that is vertically exact. $\Box $

\section{The main results}

We can state and prove the main results of the paper.

\begin{theor}
\label{thm:22} The lifted foliation ${\cal F}^{r}$ on $\nu ^{r}{\cal F}$ is
riemannian for some $r\geq 1$ iff ${\cal F}$ is a riemannian foliation.
\end{theor}

{\em Proof.} The sufficiency is well-known and it follows by Proposition \ref%
{pr0}. The necessity follows considering the submanifold inclusion $%
I_{0}^{r}(M)\subset \nu ^{r}{\cal F}$. The induced transverse riemannian
metric on the foliation induced on $I_{0}^{r}(M)$ gives a transverse
riemannian metric on $\nu {\cal F}$, thus ${\cal F}$ is riemannian. $\Box $

We say that a foliation ${\cal F}$ is transversaly almost parallelizable if
there is a ${\cal F}$-transverse vector bundle $\xi $ over $M$, such that $%
\xi \oplus \nu {\cal F}$ is transversaly parallelizable. If a foliation $%
{\cal F}$ is transversaly parallelizable, then it is a riemannian one by a
transverse metric given by a parallelization. In the case of an almost
transversaly parallelizable, any transverse riemannian metric given by a
parallelism of $\xi \oplus \nu {\cal F}$ induces a transverse riemannian
metric on $\nu {\cal F}$. Thus the following statement holds true.

\begin{cor}
\label{corpar}If the lifted foliation ${\cal F}^{r}$ on $\nu ^{r}{\cal F}$
is transversaly parallelizable or almost parallelizable for some $r\geq 1$,
then ${\cal F}$ is a riemannian foliation.
\end{cor}

The proof of Theorem \ref{thm:22} can not give any answer to the following
question: {\em when is }${\cal F}${\em \ riemannian if the foliation induced
on }$\nu ^{r}{\cal F}\backslash I_{r-1}^{r}(\nu ^{r-1}{\cal F})${\em \ is
riemannian for some }$r\geq 1$? We are going to relate this question to the
existence of a certain transverse slashed lagrangian $L^{r}$\ of order $r$,
asking that the open subset $\nu _{\ast }^{r}{\cal F\subset }\nu ^{r}{\cal F}
$ that does not contains $I_{r-1}^{r}(\nu ^{r-1}{\cal F})$. We say that a
such lagrangian $L^{r}$ is $r$-regular if its vertical hessian, according to
the induced affine bundle structure $\pi _{r-1}^{r}:\nu ^{r}{\cal F}%
\rightarrow \nu ^{r-1}{\cal F}$, is non-degenerate. In order to give an
answer to the above question, we consider below some other regularity
conditions for slashed lagrangians of order $r$.

A transverse bundle of order $r$, $\nu ^{r}{\cal F}$ can be regarded as a
fibered manifold $\pi _{r^{\prime }}^{r}:\nu ^{r}{\cal F}\rightarrow \nu
^{r^{\prime }}{\cal F}$, $(\forall )0\leq r^{\prime }<r$. We denote $\nu
^{r^{\prime },(r-r^{\prime })\ast }{\cal F}=\nu ^{r^{\prime }}{\cal F}\times
_{M}\left( \nu ^{\ast }{\cal F}\right) ^{r-r^{\prime }}$ (where $\left( \nu
^{\ast }{\cal F}\right) ^{r-r^{\prime }}=$ $\nu ^{\ast }{\cal F}\times
_{M}\cdots \times _{M}\nu ^{\ast }{\cal F}$, with the fibered product of$\
(r-r^{\prime })$-times and $\nu ^{\ast }{\cal F}$ is the transverse bundle
dual to $\nu {\cal F}$).

In particular, according to the case of trivial foliation of $M$ by points
in \cite{MP06}, $\nu ^{1.(r-1)\ast }{\cal F}=\nu ^{r-1}{\cal F}\times
_{M}\nu ^{\ast }{\cal F}$ is denoted by $\nu ^{r\ast }M$ and play the role
of the vectorial dual of the affine bundle $\nu ^{r}{\cal F}\rightarrow \nu
^{r-1}{\cal F}$.

A transverse slashed lagrangian of order $r$ is a map $L^{r}:$ $\nu ^{r}%
{\cal F}\rightarrow I\!\!R$ that is differentiable on an open subset $\nu
_{\ast }^{r}{\cal F\subset }\nu ^{r}{\cal F}$, where $\nu ^{r}{\cal F}$ $%
\backslash \nu _{\ast }^{r}{\cal F}$ contains $I_{r-1}^{r}(\nu ^{r-1}{\cal F}%
)$.

We denote by $\nu _{\ast }^{r^{\prime }}{\cal F}=\pi _{r^{\prime }}^{r}(\nu
_{\ast }^{r}{\cal F)\subset }\nu ^{r^{\prime }}{\cal F}$ and we consider the
slashed bundles $\nu _{\ast }^{\ast }{\cal F=\nu }^{\ast }{\cal F}\backslash
\{\bar{0}\}$ and $\nu _{\ast }^{r^{\prime },(r-r^{\prime })\ast }{\cal F}%
=\nu _{\ast }^{r^{\prime }}{\cal F\times }_{M}\left( \nu _{\ast }^{\ast }%
{\cal F}\right) ^{r-r^{\prime }}$, for $0\leq r^{\prime }\leq r-1$. The
elements of a fiber $\left( \nu _{\ast }^{r^{\prime },(r-r^{\prime })\ast }%
{\cal F}\right) _{;m}$ of $\nu _{\ast }^{r^{\prime },(r-r^{\prime })\ast }%
{\cal F}\rightarrow M$ are couples of higher order elements and $%
(r-r^{\prime })$ first order momenta. The usual partial derivatives of $L$
in the highest order transverse coordinates define a well-defined Legendre
map ${\cal L}^{(r-1)}:\nu _{\ast }^{r}\rightarrow \nu ^{1,(r-1)\ast }{\cal F}%
=\nu ^{r-1}{\cal F}\times _{M}\nu ^{\ast }{\cal F}$. We suppose first that $%
{\cal L}^{(r)}$ is a diffeomorphism on its image and the image is exactly $%
{\cal L}^{(r-1)}\left( \nu _{\ast }^{r}\right) =\nu _{\ast }^{1,(r-1)\ast }%
{\cal F}=\nu _{\ast }^{r-1}{\cal F}\times _{M}\nu _{\ast }^{\ast }{\cal F}$.
Then the energy $L^{(r-1)}:\nu _{\ast }^{1,(r-1)\ast }{\cal F}\rightarrow
I\!\!R$ of the dual affine hamiltonian of $L$ (see \cite{MP06} for its
classical definition and \cite{MP09B} for a coordinate description of the
whole construction in the non-foliate case). Analogous, for $0\leq j<r-1$,
we suppose, step by step, backward from $r-1$ from $0$, that there the usual
partial derivatives of $L^{(j+1)}:\nu _{\ast }^{j+1,(r-j-1)\ast }{\cal F}%
=\nu _{\ast }^{r-j-1}{\cal F}\times _{M}\left( \nu _{\ast }^{\ast }{\cal F}%
\right) ^{j+1}\rightarrow I\!\!R$ in the highest order transverse
coordinates (of order $j+1$) define a well-defined Legendre map ${\cal L}%
^{(j+1)}:\nu _{\ast }^{j+1,(r-j-1)\ast }{\cal F}=\nu _{\ast }^{j+1}{\cal F}%
\times _{M}\left( \nu _{\ast }^{\ast }{\cal F}\right) ^{r-j-1}\rightarrow
\nu ^{j,(r-j)\ast }{\cal F}=\nu ^{j}{\cal F}\times _{M}\left( \nu ^{\ast }%
{\cal F}\right) ^{r-j}$. We suppose that ${\cal L}^{(j+1)}$ is a
diffeomorphism on its image and the image is exactly ${\cal L}^{(j+1)}\left(
\nu _{\ast }^{j+1,(r-j-1)\ast }{\cal F}\right) =\nu _{\ast }^{j,(r-j)\ast }%
{\cal F}=\nu _{\ast }^{j}{\cal F}\times _{M}\left( \nu _{\ast }^{\ast }{\cal %
F}\right) ^{r-j}$. Then the energy $L^{(j)}:\nu _{\ast }^{j,(r-j)\ast }{\cal %
F}\rightarrow I\!\!R$ of the dual affine hamiltonian of $L^{(j+1)}$ can be
considered. Finally, for $j=0$, we obtain a transverse lagrangian $%
L^{(0)}:\nu _{\ast }^{0,r\ast }{\cal F}=\left( \nu _{\ast }^{\ast }{\cal F}%
\right) ^{r}\rightarrow I\!\!R$ as the energy of the dual affine hamiltonian 
$L^{(1)}$ and we suppose that ${\cal L}^{(1)}$ $:\nu _{\ast }^{1,(r-1)\ast }%
{\cal F}=\nu _{\ast }{\cal F}\times _{M}\left( \nu _{\ast }^{\ast }{\cal F}%
\right) ^{r-1}\rightarrow \nu _{\ast }^{0,r\ast }{\cal F}=\left( \nu _{\ast
}^{\ast }{\cal F}\right) ^{r}\subset \nu ^{0,r\ast }{\cal F}=\left( \nu
^{\ast }{\cal F}\right) ^{r}$ is a diffeomorphism. It follows a
diffeomorphism ${\cal L=L}^{(1)}\circ \cdots \circ {\cal L}^{(r)}:\nu _{\ast
}^{r}\rightarrow \left( \nu _{\ast }^{\ast }{\cal F}\right) ^{r}$ and a
transverse lagrangian $L^{(0)}:\left( \nu _{\ast }^{\ast }{\cal F}\right)
^{r}\rightarrow I\!\!R$. The canonical diagonal inclusion $\nu ^{\ast }{\cal %
F}\rightarrow \left( \nu ^{\ast }{\cal F}\right) ^{r}$ sends $\nu _{\ast
}^{\ast }{\cal F}\rightarrow \left( \nu _{\ast }^{\ast }{\cal F}\right) ^{r}$%
. We suppose that the restriction of $L^{(0)}$ to the diagonal is a
positively admissible lagrangian on $\nu ^{\ast }{\cal F}$, in fact a
transverse hamiltonian $H:\nu ^{\ast }{\cal F}\rightarrow I\!\!R$. If the
given transverse lagrangian $L^{r}:\nu ^{r}{\cal F}\rightarrow I\!\!R$
fulfills all the above conditions, we say that $L$ itself is a {\em %
positively admissible lagrangian} (of order $r$) and $H$ is its{\em \
diagonal} {\em hamiltonian}.

The vertical bundle $V$ of the affine bundle $\pi _{r-1}^{r}:\nu ^{r}{\cal F}%
\rightarrow \nu ^{r-1}{\cal F}$ has the form $\left( \pi _{0}^{r}\right)
^{\ast }\nu {\cal F}\rightarrow \nu ^{r}{\cal F}$. We say that a transverse
metric $g$ on $\nu ^{r}{\cal F}$ is {\em vertically exact} if there is a
positively admissible lagrangian of order $r$, $L:\nu ^{r}{\cal F}%
\rightarrow I\!\!R$ such that the restriction of $g$ to $V$ is the same as
the vertical hessian ${\rm Hess}\ L$. We say in this case that the
riemannian foliation ${\cal F}^{r}$ is {\em vertically exact}. These
definitions can be easily adapted for the case of slashed bundles $\nu
_{\ast }^{r}{\cal F}$.

The main technical tool to prove the necessity of Theorems \ref{thm:22-1}
and \ref{thm:21} below is the following result proved in \cite[Proposition
2.2]{PMCR}.

\begin{prop}
\label{prop:23} Let $p_{1}:E_{1}\rightarrow M$ and $p_{2}:E_{2}\rightarrow M$
be foliated vector bundles over a foliated manifold $(M,{\cal F})$ and $%
q_{2}:E_{2\ast }\rightarrow M$ be the slashed bundle. If there are a
positively admissible lagrangian $L:E_{2}\rightarrow I\!\!R$ and a metric $b$
on the pull back bundle $q_{2}^{\ast }E_{1}\rightarrow E_{2\ast }$, foliated
with respect to ${\cal F}_{E_{2\ast }}$, then there is a foliated metric on $%
E_{1}$, with respect to ${\cal F}$.
\end{prop}

We can now state and prove the following Theorems.

\begin{theor}
\label{thm:22-1} Let ${\cal F}$ be a foliation on a manifold $M$ and ${\cal F%
}_{0}^{r}$ be the lifted foliation in a suitable slashed bundle $\nu _{\ast
}^{r}{\cal F}$ of the $r$-normal bundle $\nu ^{r}{\cal F}$. Then ${\cal F}%
_{0}^{r}$ is riemannian and vertically exact for some $r\geq 1$ iff ${\cal F}
$ is riemannian.
\end{theor}

In particular, it follows that any transverse metric $g$ on $\nu F$ gives
rise to a canonical lagrangian on ${\cal \nu }^{r}{\cal F}$, coming from the
vertical part of the vertically exact transverse riemannian metric on $\nu 
{\cal F}^{r}$. So, it is natural to ask that only the existence of a
lagrangian on ${\cal \nu }^{r}{\cal F}$ guaranties that ${\cal F}$ is
riemannian. One have a positive answer, as follows.

\begin{theor}
\label{thm:21} If $(M,{\cal F})$ is a foliated manifold, then there is a
positively admissible lagrangian on ${\cal \nu }^{r}{\cal F}$ for some $%
r\geq 1$ iff the foliation ${\cal F}$ is riemannian.
\end{theor}

{\em Proof} ({\em of Theorems \ref{thm:22-1} and \ref{thm:21}}). The
sufficiency for both Theorems follow by Proposition \ref{pr0}. The necessity
for Theorem \ref{thm:21} follows using Proposition \ref{prop:23} with $%
E_{1}=E_{2}=\nu ^{\ast }{\cal F}$ and $H$ the diagonal hamiltonian. Finally,
the necessity for Theorem \ref{thm:22-1} follows thanks to Theorem \ref%
{thm:21} for the positively admissible lagrangian on ${\cal \nu }^{r}{\cal F}
$, given by the condition that the riemannian metric on ${\cal \nu }^{r}%
{\cal F}$ is vertically exact. $\Box $

Finally, as in \cite{PMCR}, the following question arises: {\em can we drop
in Theorem \ref{thm:22-1} the condition that }${\cal F}_{0}^{r}${\em \ is
vertically exact?}


\begin{thebibliography}{99}
\bibitem{Bu} Bucataru I., {\em Canonical semisprays for higher order
Lagrange spaces}, C. R. Acad. Sci. Paris, Ser. I 345 (2007) 269--272.

\bibitem{CSC1} Crampin M., Sarlet W., Cantrijn F., {\em Higher Order
differential equations and higher order lagrangian mechanics}, Math. Proc.
Camb. Phil. Soc., 86 (1986), 565--587.

\bibitem{Mi} R. Miron, {\em The Geometry of Higher-Order Lagrange Spaces.
Applications to Mechanics and Physics}, FTPH no. 82, Kluwer Academic
Publisher, 1997.

\bibitem{JoWo} J\'{o}zefowicz M., Wolak R., {\em Finsler foliations of
compact manifolds are Riemannian}, Differential Geometry and its
Applications, 26 (2) (2008) 224--226.

\bibitem{MiMo} Miernowski A., Mozgawa W., {\em Lift of the Finsler foliation
to its normal bundle}, Differential Geometry and its Applications, 24 (2006)
209--214.

\bibitem{Mo} Molino P., {\em Riemannian foliations}, Progress in
Mathematics, Vol. 73, Birh\"{a}user, Boston, 1988.

\bibitem{Mo1} Morimoto A., {\em Prolongations of G-structure to tangent
bundles of higher order}, Nagoya Math. J., 38 (1970) 153-179.

\bibitem{Mo2} Morimoto, A., {\em Liftings of tensor fields and connections
to tangent bundles of higher order}, Nagoya Math. J., 40 (1970) 99-120.

\bibitem{MP06} Popescu P., Popescu M., {\em Affine Hamiltonians in higher
order geometry}, International Journal of Theoretical Physics, 46 (10)
(2007) 2531-2549.

\bibitem{MP09F} Popescu P., Popescu M., {\em Lagrangians adapted to
submersions and foliations}, Differential Geometry and its Applications, 27
(2) (2009) 171-178.

\bibitem{MP09B} Popescu M., Popescu P., {\em Lagrangians and higher order
tangent spaces}, Balkan Journal of Geometry and its Applications, 15 (1)
(2010) 142-148.

\bibitem{PMCR} Popescu P., Popescu M., {\em Foliated vector bundles and
Riemannian foliations}, C. R. Acad. Sci. Paris, Ser. I 349 (2011) 445--449.

\bibitem{Ta01} Tarquini C., {\em Feuilletages de type fini compact}, C. R.
Acad. Sci. Paris, Ser. I, 339 (2004) 209-214.

\bibitem{Va} Vaisman I., {\em Hamiltonian structures on foliations}, J.
Math. Phys., 43, 10 (2002), 4966-4977.

\bibitem{Wo1} Wolak R.A., {\em On transverse structures of foliations},
Proceedings of the 13th winter school on abstract analysis (Srn\'{\i},
1985). Rend. Circ. Mat. Palermo (2) Suppl. No. 9 (1985), 227--243.

\bibitem{Wo2} Wolak R.A., {\em Leaves of foliations with a transverse
geometric structure of finite type}, Publ. Mat. 33 (1989), no. 1, 153--162

\bibitem{Wo3} Wolak R.A., {\em Foliated and associated geometric structures
on foliated manifolds}, Ann. Fac. Sci. Toulouse, V. S\'{e}r., Math., 10 (3)
(1989) 337-360.
\end{thebibliography}
\end{document}